\documentclass[conference]{IEEEtran}
\usepackage{cite}
\usepackage{amsmath,amssymb,amsfonts}
\usepackage{algorithmic}
\usepackage{graphicx}
\usepackage{textcomp}
\usepackage{xcolor}
\usepackage{multirow}

\usepackage{bm}
\def\BibTeX{{\rm B\kern-.05em{\sc i\kern-.025em b}\kern-.08em
    T\kern-.1667em\lower.7ex\hbox{E}\kern-.125emX}}
\begin{document}

\title{Efficient, arbitrarily high precision hardware logarithmic arithmetic for linear algebra\\
}

\author{\IEEEauthorblockN{Jeff Johnson}
\IEEEauthorblockA{\textit{Facebook AI Research} \\
Los Angeles, CA, United States \\
jhj@fb.com}
}


\maketitle

\begin{abstract}
  The logarithmic number system (LNS) is arguably not broadly used due to exponential circuit overheads for summation tables relative to arithmetic precision. Methods to reduce this overhead have been proposed, yet still yield designs with high chip area and power requirements. Use remains limited to lower precision or high multiply/add ratio cases, while much of linear algebra (near 1:1 multiply/add ratio) does not qualify.

We present a dual-base approximate logarithmic arithmetic comparable to floating point in use, yet unlike LNS it is easily fully pipelined, extendable to arbitrary precision with $\mathcal{O}(n^2)$ overhead, and energy efficient at a 1:1 multiply/add ratio. Compared to float32 or float64 vector inner product with FMA, our design is respectively 2.3$\times$ and 4.6$\times$ more energy efficient in 7 nm CMOS. It depends on exp and log evaluation 5.4$\times$ and 3.2$\times$ more energy efficient, at 0.23$\times$ and 0.37$\times$ the chip area for equivalent accuracy versus standard hyperbolic CORDIC using shift-and-add and approximated ODE integration in the style of Revol and Yakoubsohn. This technique is a novel alternative for low power, high precision hardened linear algebra in computer vision, graphics and machine learning applications.

\end{abstract}


\section{Introduction}\label{sec:intro}

Energy efficiency is typically the most important challenge in advanced CMOS technology nodes. With the \textit{dark silicon} problem~\cite{Venkatesh:2010:CCR:1736020.1736044}, the vast majority of a large scale design is clock or power gated at any given point in time. Chip area becomes exponentially more available relative to power consumption, preferring ``a new class of architectural techniques that 'spend' area to 'buy' energy efficiency''~\cite{taylor2012dark}. Memory architecture is often the most important concern, with 170-6400$\times$ greater DRAM access energy versus arithmetic at 45 nm~\cite{horowitz20141}.
This changes with the rise of machine learning, as heavily employed linear algebra primitives such as matrix/matrix product offer substantial local reuse of data by algorithmic tiling~\cite{Wolfe:1989:MIS:76263.76337}: $\mathcal{O}(n^3)$ arithmetic operations versus $\mathcal{O}(n^2)$ DRAM accesses.
This is a reason for the rise of dedicated neural network accelerators, as memory overheads can be substantially amortized over many arithmetic operations in a fixed function design, making arithmetic efficiency matter again.

Many hardware efforts for linear algebra and machine learning tend towards low precision implementations, but here we concern ourselves with the opposite: enabling (arbitrarily) high precision yet energy efficient substitutes for floating point or long word length fixed point arithmetic. There are a variety of ML, computer vision and other algorithms where accelerators cannot easily apply precision reduction, such as hyperbolic embedding generation~\cite{nickel2017poincare} or structure from motion via matrix factorization~\cite{Tomasi:1992:SMI:144398.144403}, yet provide high local data reuse potential.

The logarithmic number system (LNS)~\cite{swartzlander1975sign} can provide energy efficiency by eliminating hardware multipliers and dividers, yet maintains significant computational overhead with \textit{Gaussian logarithm} functions needed for addition and subtraction. While reduced precision cases can limit themselves to relatively small LUTs/ROMs, high precision LNS require massive ROMs, linear interpolators and substantial MUXes. Pipelining is difficult, requiring resource duplication or handling variable latency corner cases as seen in~\cite{coleman2008processor}. The ROMs are also exponential in LNS word size, so become impractical beyond a float32 equivalent. Chen et al.~\cite{Chen2000} provide an alternative fully pipelined LNS add/sub with ROM size a $\mathcal{O}(n^3)$ function of LNS word size, extended to float64 equivalent precision. However, in their words, ``[our] design of [a] large word-length LNS processor becomes impractical since the hardware cost and the pipeline latency of the proposed LNS unit are much larger.'' Their float64 equivalent requires 471 Kbits ROM and at least 22,479 full adder (FA) cells, and 53.5 Kbits ROM and 5,550 FA cells for float32, versus a traditional LNS implementation they cite~\cite{Lewis1994} with 91 Kbits of ROM and only 586 FA cells.

While there are energy benefits with LNS~\cite{Popoff:2016:HLN:2971808.2972130}, we believe a better bargain can be had. Our main contribution is a trivially pipelined logarithmic arithmetic extendable to arbitrary precision, using no LUTs/ROMs, and a $\mathcal{O}(n^2)$ precision to FA cell dependency. Unlike LNS, it is substantially more energy efficient than floating point at a 1:1 multiply/add ratio for linear algebra use cases. It is approximate in ways that an accurately designed LNS is not, though with parameters for tuning accuracy to match LNS as needed. It is based on the \textit{ELMA} technique~\cite{jhj}, extended to arbitrary precision with an energy efficient implementation of exp/log using restoring shift-and-add~\cite{Muller:1985:DBC:4135.4987} and an ordinary differential equation integration step from Revol and Yakoubsohn\cite{Revol1999AcceleratedSA} but with approximate multipliers and dividers. It is tailored for vector inner product, a foundation of much of linear algebra, but remains a general purpose arithmetic.
We will first describe our hardware exp/log implementations, then detail how they are used as a foundation for our arithmetic, and provide an accuracy analysis. Finally, hardware synthesis results are presented and compared with floating point.

\section{Notes on hardware synthesis}

All designs considered in this paper are on a commercially available 7 nm CMOS technology constrained to only SVT cells. They are generated using Mentor Catapult high level synthesis (HLS), biased towards min latency rather than min area, with ICG (clock gate) insertion where appropriate.
Area is reported via Synopsys Design Compiler, and power/energy is from Synopsys PrimeTime PX from realistic switching activity. Energy accounts for combinational, register, clock tree and leakage power, normalized with respect to module throughput in cycles, so this is a per-operation energy. We consider pipelining acceptable for arithmetic problems in linear algebra with sufficient regularity such as matrix multiplication (Section~\ref{sec:flma}), reducing the need for purely combinational latency reduction.
Power analysis is at the TT@25C corner at nominal voltage. Design clocks from 250-750 MHz were considered, with 375 MHz chosen for reporting, being close to minimum energy for many of the designs. Changing frequency does change pipeline depth and required register/clock tree power, as well as choice of inferred adder or other designs needed to meet timing closure by synthesis.

\section{exp/log evaluation}

Our arithmetic requires efficient hardware implementation of exponential $b^x$ and logarithm $\log_b(x)$ for a base $b$, which are useful in their own right.
Typical algorithms are power series evaluation, polynomial approximation/table-based methods, and shift-and-add methods such as hyperbolic CORDIC~\cite{Volder:1959:CCT:1457838.1457886} or the simpler method by De Lugish~\cite{DeLugish:1970:CAA:905445}.
Hardware implementations have been considered for CORDIC~\cite{Duprat:1993:CAN:626493.626709}, ROM/table-based implementations~\cite{Schulte:1994:HDE:626509.626954}, approximation using shift-and-add~\cite{Abed:2003:CVI:951846.951880} with the Mitchell logarithm approximation~\cite{5219391}, and digit recurrence/shift-and-add~\cite{Pineiro:2005:HLS:1050556.1050588}. CORDIC requires three state variables and additions per iteration, plus a final multiplication by a scaling factor.
BKM~\cite{295857} avoids the CORDIC scaling factor but introduces complexity in the iteration step.

Much of the hardware elementary function literature is concerned with latency reduction rather than energy optimization. Variants of these algorithms such as high radix formulations~\cite{Baker:1975:PMA:1309297.1309908}\cite{Ercegovac93}~\cite{Pineiro:2005:HLS:1050556.1050588} or parallel iterations~\cite{Duprat:1993:CAN:626493.626709} increase switching activity via additional active area, iteration complexity, or adding sizable MUXes in the radix case. In lieu of decreasing combinational delay via parallelism, pipelining is a worthwhile strategy to reduce energy-delay product~\cite{Gonzalez1996}, but only with high pipeline utilization and where register power increases are not substantial. Ripple-carry adders, the simplest and most energy efficient adders, remain useful in the pipelined regime, and variants like variable block adders improve latency for minimal additional energy~\cite{Vratonjic2005}. Fully parallel adders like carry-save can improve on both latency and switching activity for elementary functions~\cite{Pineiro:2005:HLS:1050556.1050588}, but only where the redundant number system can be maintained with low computational overhead. For example, in shift-and-add style algorithms, adding a shifted version of a carry-save number to itself requires twice the number of adder cells as a simple ripple-carry adder (one to add each of the shifted components), resulting in near double the energy. Eliminating registers via combinational multicycle paths (MCPs) is another strategy, but as the design is no longer pipelined, throughput will suffer, requiring an introduction of more functional units or accepting the decrease in throughput.
There is then a tradeoff between clock frequency, combinational latency reduction, pipelining for timing closure, MCP introduction, and functional unit duplication versus energy per operation.




\section{$e^x$ shift-and-add with integration}\label{sec:exp}

We consider De Lugish-style restoring shift-and-add, which will provide ways to reduce power or recover precision with fewer iterations (Section~\ref{sec:integrate} and~\ref{sec:approx-mul}).
The procedure for exponentials $y=b^x$ is described in Muller~\cite{Muller:1985:DBC:4135.4987} as:
\begin{align*}
  L_0 &= x \\
  L_{n+1} &= L_n - \log_b(1 + d_n 2^{-n}) \\
  E_0 &= 1 \\
  E_{n+1} &= E_n(1 + d_n 2^{-n}) \\
  d_n &= \begin{cases}
    1 & \text{if } L_n \geq \log_b(1 + 2^{-n}) \\
    0 & \text{otherwise}
  \end{cases} \\
  y &= E_I ~(\text{at desired iteration $I$})
\end{align*}
The acceptable range of $x$ is $[0, \sum^{\infty}_{n=0} \log_b(1 + 2^{-n}))$, or $[0, 1.56\dots)$ for $b=e$ (Euler's number). Range reduction techniques considered in~\cite{Muller:1985:DBC:4135.4987} can be used to reduce arbitrary $x$ to this range. This paper will only consider $b=e$ and limiting $x$ as fixed point, $x \in [0, \ln(2))$, restrictions discussed in Sections~\ref{sec:integrate} and~\ref{sec:dual}.

      We must consider rounding error and precision of $x$, $L_n$ and $E_n$. Our range-limited $x$ can be specified purely as a fixed point fraction with $x_{\text{bits}}$ fractional bits. The iteration $n=1$ is skipped as $x \in [0, \ln(2))$. All subsequent $L_n$ are $<1$ and can be similarly represented as a fixed point fraction.
          These $L_n$ will use $\ell$ fractional bits ($\ell \geq x_{\text{bits}}$) with correctly rounded representations of $\ln(1 + d_n 2^{-n})$. 
          $E_n \in [1, 2)$ is the case in our restricted domain, which is maintained as a fixed point fraction with an implicit, omitted leading integer 1. Multiplication by $2^{-n}$ is a shift by $n$ bits, so we append this leading 1 to the post-shifted $E_n$ before addition. We use $p$ fractional bits to represent $E_n$. At the final $I$-th iteration, $E_I$ is rounded to $y_{\text{bits}} \leq p$ bits for the output $y$. Ignoring rounding error, the relative error of the algorithm is $|e^x - E_I| / e^x = 2^{-I+1}$ at iteration $I$, so for 23 fractional bits, $I=24$ is desired.

            All adders need not be of size $\ell$ or $p$, either. $L_n$ reduces in magnitude at each step; with $y \in [1, 2)$, $L_n$ only needs the $\ell - \max(0, n-2)$ LSB fractional bits.
$E_n$ has a related bit size progression $0, 0, 1, 3, 5, 9, 14, \dots$, as $E_n$ is $\max(p, \text{size}(E_{n-1}) + n)$ fractional bits, except starting at $n \geq 3$ we are off by 1 ($d_1, d_2, d_3$ cannot all be 1, as $\sum_{i=1}^{3} \ln(1 + 2^{-i}) > \ln(2)$).
While $L_n$ successively reduces in precision from $\ell$, we limit $E_n$ to $p$ fractional bits via truncation (bits shifted beyond position $p$ are ignored). As with~\cite{Pineiro:2005:HLS:1050556.1050588}, we can deal with truncation error by setting $p = y_{\text{bits}} + (\lceil \log_2(I) \rceil + 1)$, using extra bits as guard bits.

$L_n$ requires an adder and MUX ($L_n = L_{n-1} - \ln(1 + 2^{-n+1})$ if $d_n = 1$, or $L_{n} = L_{n-1}$ if $d_n = 0$). The constants $\ln(1 + 2^{-n})$ are hard-wired into adders when iterations are fully unrolled (a separate adder for each iteration).
The $E_n$ do not use a full adder of size $p$ in the general case; only shifted bits that overlap with previously available bits need full adder cells. The $L_n$ can also be performed first, with $d_n$ stored in flops (to reduce glitches) for data gating $E_n$ additions, reducing switching activity at the expense of higher latency, as 25\% of the $d_i$ on average will remain zero across iterations.

One can use redundant number systems for $L_n$ and $E_n$ and avoid full evaluation of the $L_n$ comparator~\cite{Muller:1985:DBC:4135.4987}, but $E_n$ is problematic. In the non-redundant case, only a subset of the shifted $E_n$ require a full adder, and the remainder only a half adder. With a carry-save representation for $E_n$, two full adders are required for the entire length of the word, one to add each portion of the shifted carry-save representation. While the carry-save addition is constant latency, it requires more full adder cells. In our evaluation carry-save for $E_n$ prohibitively increases power over synthesis-inferred adder choice. At high clock frequencies (low latency) this tradeoff is acceptable, but low power designs will generally avoid this regime.

\begin{figure}
    \vspace{-0.6cm}
  \begin{center}
    \includegraphics[width=1.0\linewidth]{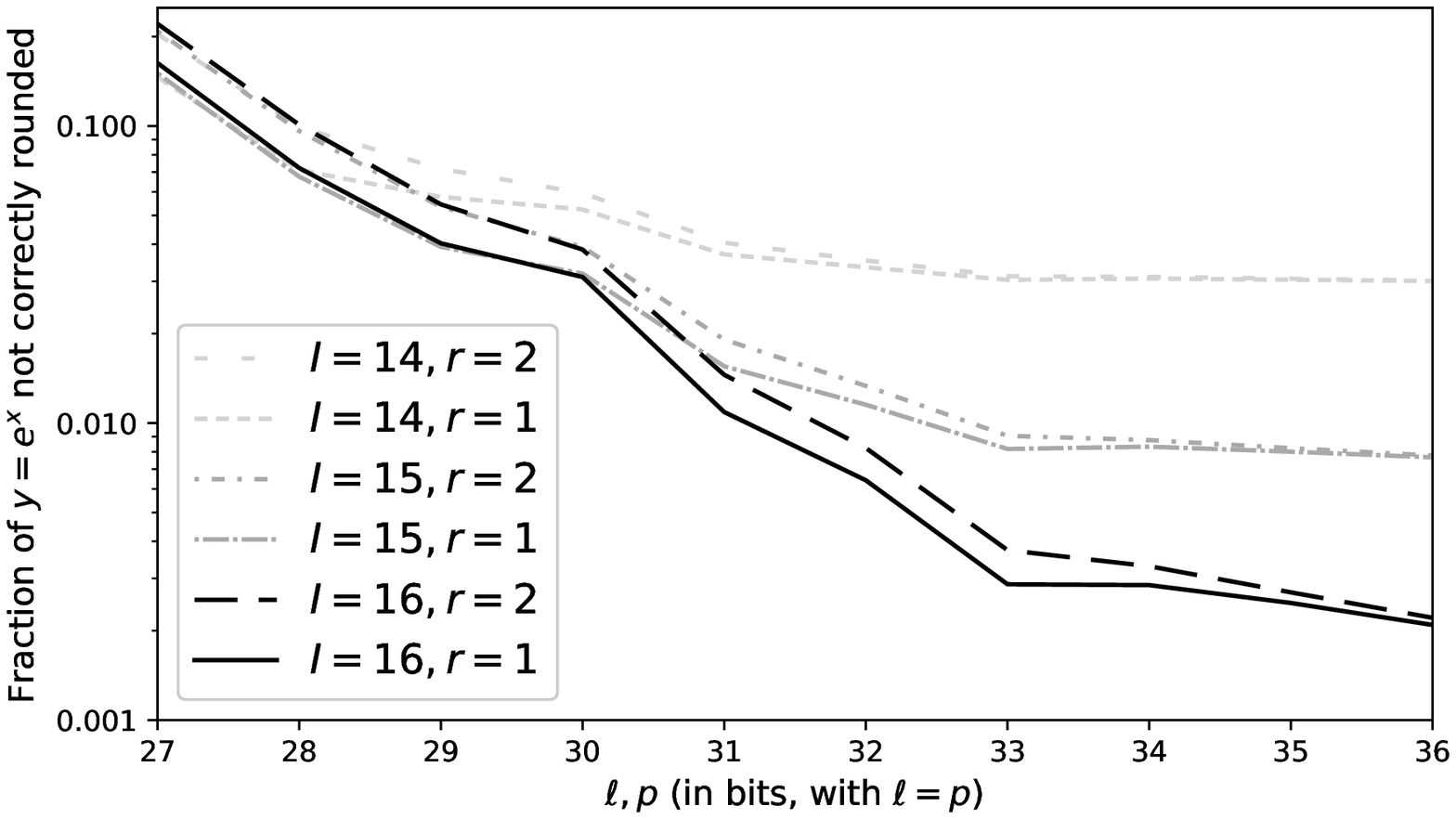}
  \end{center}
  \vspace{-0.5cm}
  \caption{\label{fig:exp_accuracy}
    Our $e^x$ accuracy at $x_{\text{bits}} = y_{\text{bits}} = 23$ relative to $I, \ell, p, r$. All configurations have $\leq1$ ulp error. }
\end{figure}
\subsection{Euler method integration}\label{sec:integrate}

This algorithm is simple but has high latency from the sequential dependency of many adders, high $\ell, p$ and iterations $I$ for accurate results.
For significant latency and energy reduction, Revol and Yakoubsohn~\cite{Revol1999AcceleratedSA} show that about half the iterations can be omitted by treating the problem as an ordinary differential equation with a single numerical integration step.
$e^x$ satisfies the ODE $y' = f(x, y) = y$ where $y(0) = 1$. They consider in software both an explicit Euler method and 4th order Runge-Kutta (RK4). RK4 involves several multipliers and is not a good energy tradeoff to avoid more iterations. The explicit Euler method step has a single multiplication:
\begin{align*}
  y &= E_I &~\text{no integration step} \Rightarrow \\
  y &= E_I + L_I E_I &~\text{explicit Euler method step}
\end{align*}
at the $I$-th terminal iteration, with residual $L_I$ used as the step size. They give a formula for $I$ at a desired accuracy of $\epsilon$ of $L_I \leq \sqrt{2\epsilon / e^{1.56}}$, ignoring truncation error. Note that $L_{n+1} < 2 \ln(1 + 2^{-n}) \approx 2^{-n+1}$. Thus, for single-precision $\epsilon = 2^{-24}$, we need $L_I \leq 2^{-12.63\dots}$, or $I \geq 14$. Double-precision $\epsilon=2^{-53}$ has $I \geq 29$, and quad precision $\epsilon=2^{-113}$ has $I \geq 59$. Implementation of $2^x$ from this requires pre-multiplication of $x$ by a fixed point rounding of $\ln(2)$, a significant energy overhead.

\subsection{Integration via approximate multiplication}\label{sec:approx-mul}

We have $L_I \in [0, 1)$, $E_I \in [1, 2)$ when $x \in [0, \ln(2))$. The Euler method step multiplication $L_I E_I$ would be a massive $\ell \times (1 + p)$ bits, with the $1+$ for the leading integer 1 bit of $E_I$. Let $L^f_I$ and $E^f_I$ denote fractional portions of $L_I$ and $E_I$. The step can then be expressed as:
      $$y = E_I + (0 + L^f_I)(1 + E^f_I) = E_I + L^f_I + L^f_I E^f_I$$
$L^f_I E^f_I$ now solely involves fractional bits, of which we only care about $y_{\text{bits}}$ to $p$ MSBs produced.
$L^f_I$ has $\max(I-2, 0)$ zero MSBs, so there are $\max(I-2, 0)$ ignorable zero MSBs in the resulting product, yielding a $p \times (\ell - \max(I-2, 0))$ multiplier, still an exact step calculation. Assuming $I > 2$ and given these zero MSBs, we only need $(p - I + 2)$ MSBs of the result, so we truncate both $L^f_I$ and $E^f_I$ to limit the result to this size (truncation ignores carries from the multiplication of the truncated LSBs). We do this symmetrically, and since usually $p > y_{\text{bits}}$, we take $\ell - (I-2) - r$ fractional MSBs from $E^f_I$, with an option to remove another $r$ bits, $0 \leq r \leq 4$. This may not produce enough bits to align properly with $p$, so we append zeros to the LSBs as needed to match the size of $p$.
      For example, at $x_{\text{bits}}=23$, $y_{\text{bits}}=24$, $\ell = p = 28$, $I=14$, $r=2$, we have a 14 $\times$ 14 multiplier, of which we only need 16 MSBs (based on alignment with $E^f_I$), and the ultimate carry from the 12 LSBs. One can consider other approximate multipliers~\cite{Jiang2016}, but truncation seems to work well and provides a significant reduction in energy.


\begin{table}
{
  \caption{\label{tab:exp-log-syn}
    Fully pipelined exp/log, $x_{\text{bits}} = y_{\text{bits}} = 23$ synthesis results
}}
\centering \begin{tabular}{|l|l|l|l|l|}
\hline
\bm{$e^x$} & \textbf{$\textbf{(0.5, 1]}$ ulp err} & \textbf{Cycles} & \textbf{Area} $\mu \text{m}^2$ & \textbf{Energy} pJ \\
\hline
CORDIC & 9.98\% & 4 & 1738 & 2.749\\
\textbf{Ours} & \textbf{9.90\%} & \textbf{2} & \textbf{407.2} (0.23$\times$) & \textbf{0.512} (0.19$\times$)\\
\hline
\textbf{\bm{$\ln(x)$}} & & & &  \\
\hline
CORDIC & 14.4\% & 4 & 2084 & 3.573\\
\textbf{Ours} & \textbf{14.8\%} & \textbf{4} & \textbf{769.4} (0.37$\times$) & \textbf{1.107} (0.31$\times$)\\
\hline
\end{tabular}
\vspace{-0.5cm}
\end{table}

\subsection{Error analysis and synthesis results}

The table maker's dilemma is unavoidable for transcendental functions~\cite{Lefevre1998}. For $x_{\text{bits}} = y_{\text{bits}} = 23$, we need to evaluate $e^x$ to at least 42 bits to provide correctly rounded results for fixed point $x \in[0, \ln(2))$. In lieu of exact evaluation, we demand function monotonicity, $\leq1$ ulp error, and consider the occurrence of incorrectly rounded results ($>0.5$ ulp error). Figure~\ref{fig:exp_accuracy} considers error in this regime with a sweep of $I, \ell$, $p$ and $r$, with $\ell = p$. $\ell, p < 27$ has maximum $>1$ ulp error.

Table~\ref{tab:exp-log-syn} shows fully-pipelined (iterations unrolled), near iso-accuracy synthesis results for our method ($I=14$, $\ell = p = 28$, $r=2$) and standard hyperbolic CORDIC (28 iterations and 29 fractional bit variables). All implementations have $\leq 1$ ulp error, with the fraction at $(0.5, 1]$ error shown shown. We are 5.4$\times$ more energy efficient, 0.23$\times$ the area, and half the latency in cycles; as discussed earlier, most CORDIC modifications reduce latency at the expense of increased energy.

\section{$\ln(x)$ shift-and-add with integration}

$y=\ln(x)$ is similar to $e^x$ with roles of $E_n$ and $L_n$ reversed, with division for the integration~\cite{Revol1999AcceleratedSA}:
\begin{align*}
  E_0 &= 1 \\
  E_{n+1} &= E_n(1 + d_n 2^{-n}) \\
  L_0 &= 0 \\
  L_{n+1} &= L_n + d_n \ln(1 + 2^{-n}) \\
  d_n &= \begin{cases}
    1 & \text{if } E_n(1 + 2^{-n}) \leq x \\
    0 & \text{otherwise}
  \end{cases} \\
  y &= L_I + (x - E_I) / E_I ~(\text{Euler method step})
\end{align*}
We restrict ourselves to $x \in [1, 2)$. The error of $E_I \approx \ln(x)$ is $\leq 2^{-I+1}$, with the target number of iterations $I$ (ignoring truncation error) for error $\epsilon$ given when $(x - E_I) \leq \sqrt{2\epsilon}$. For single precision $\epsilon = 2^{-24}$, $I=13$, and double precision $\epsilon = 2^{-53}$, $I=27$, and $\epsilon = 2^{-113}$ is $I=57$.
  Prior discussion concerning the $E_n$ and $L_n$ sequences and data gating with $d_n$ carry over to this algorithm. It is also the case that the running sum $L_n$ is not needed until the very end, so a carry-save adder postponing full evaluation of carries is appropriate.
It is possible to use a redundant number system for $E_I$ and avoid full evaluation of the comparison~\cite{Muller:1985:DBC:4135.4987}, but the required shift with add increases switching activity significantly.

\subsection{Integration via approximate division}

We approximate the integration division by truncating the dividend and divisor. The dividend $(x - E_I) \in [0, 1)$ has at least $\max(0, I - 3)$ zero fractional MSBs, and the divisor $E_I \in [1, 2)$, so the result is a fraction that we must align with the $\ell$ bits in $L_I$ for the sum.
  We skip known zero MSBs, and some number $r$ of the LSBs of the dividend. For the divisor $E_I$, we need not use the entire fractional portion but choose only some number of fractional bits $s$. We then have a $(p - \max(0, I-3) - r)$ by $1+s$ fixed point divider ($1+$ is for the leading integer 1 of $E_I$). $r = 3, s = 9$ is reasonable in our experiments. This is higher area and latency than the truncated multiplier (we only evaluated truncated division with digit recurrence), but the increase in resources of log versus exp is acceptable for linear algebra use cases (Section~\ref{sec:syn}).

  \begin{figure}
    \vspace{-0.6cm}
  \begin{center}
    \includegraphics[width=1.0\linewidth]{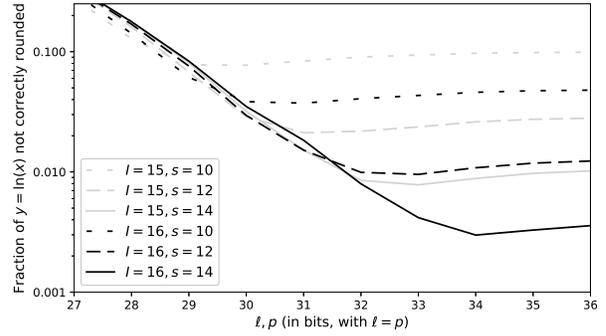}
  \end{center}
  \vspace{-0.5cm}
  \caption{\label{fig:log_accuracy}
    Our $\ln(x)$ accuracy at $x_{\text{bits}} = y_{\text{bits}} = 23$ relative to $I, \ell, p, s$ with $r=3$. All configurations have $\leq 1$ ulp error.
  }
\end{figure}
\subsection{Error analysis and synthesis results}

As before, we only consider monotonic implementations with $\leq 1$ ulp error, and consider the frequency of incorrectly rounded results. Figure~\ref{fig:log_accuracy} shows such error occurrence versus a sweep of $I$, $\ell$, $p$, $s$, with $\ell = p$. $s$ has a larger accuracy effect than $r$, and $r=3$ yields reasonable results, so all are constrained to $r=3$. Table~\ref{tab:exp-log-syn} shows near iso-accuracy synthesis results for our method ($I=15$, $\ell = p = 28$, $r=3$, $s=9$) and standard hyperbolic CORDIC (28 iterations and 30 fractional bit variables). Our implementation is 3.2$\times$ more energy efficient at 0.37$\times$ area versus CORDIC, with much of the latency and energy coming from the truncated divider. The higher resource consumption of log over exp CORDIC is from the initialization of the $X$ and $Y$ CORDIC variables to $x+1$, $x-1$ rather than 1, with propagation of inferred required adder length by HLS throughout when synthesizing the design.

\section{Approximate logarithmic arithmetic}\label{sec:approx}

We show how the preceding designs are used to build an arbitrarily high precision logarithmic arithmetic with some (tunably) approximate aspects.

\subsection{LNS arithmetic}

The sign/magnitude logarithmic number system (LNS)~\cite{swartzlander1975sign} represents values $x \in \mathbb{R}$ as a rounded fixed point representation to some number of integer and fractional bits of $\log_b(|x|)$, plus a sign and zero flag.
The base $b$ is typically 2. We refer to $x' = \{\pm \log_b(|x|)~\text{or}~0\}$ as a representation of $x$ in the \textit{log domain}. We refer to rounding and encoding $x$ as integer, fixed or floating point as a \textit{linear domain} representation, though note that floating point itself is a combination of log and linear representations for the exponent and significand.

The benefit of LNS is simplifying multiplication, division and power/root.
For log domain $x', y'$, multiplication or division of the corresponding linear domain $x$ and $y$ is $x' \pm y'$, $n$-th power of $x$ is $nx'$ and $n$-th root of $x$ is $x'/n$, with sign, zero (and infinity/NaN flags if desired) handled in the obvious manner.
Addition and subtraction, on the other hand, require \textit{Gaussian logarithm} computation. For linear domain $x, y$, log domain add/sub of the corresponding $x', y'$ is:
\begin{align*}
\log_b(|x| + |y|) &= x' + \log_b(1 + b^{r}) \\
\log_b(|x| - |y|) &= x' + \log_b(|1 - b^{r}|)
\end{align*}
where $r = y' - x'$. Without loss of generality, we restrict $y' \leq x'$, so we only consider $r \leq 0$.
These functions are usually implemented with ROM/LUT tables (possibly with interpolation) rather than direct function evaluation, ideally realized to $\leq$ 0.5 log domain ulp relative error. The subtraction function has a singularity at $r = 0$, corresponding to exact cancellation $y - x = 0$, with the region $r \in [-\epsilon, 0)$ very near the singularity corresponding to near-exact cancellation. Realizing this critical region to 0.5 log ulp error without massive ROMs (241 Kbits in~\cite{coleman:7061396}) is a motivation for \textit{subtraction co-transformation} to avoid the singularity, which can reduce the requirement to at least 65 Kbits~\cite{Popoff:2016:HLN:2971808.2972130}. Some designs are proposed as being ROM-less~\cite{Ismail2011}, but in practice the switching power and leakage of the tables' combinational cells would still be huge. Interpolation with reduced table sizes can also be used, but the formulation in~\cite{Taylor1983} only considers log addition without the singularity.  An ultimate limit on the technique not far above float32 equivalent is still faced, as accurate versions of these tables scale exponentially with word precision~\cite{Chen2000}.

Pipelined LNS add/sub is another concern. As mentioned in Section~\ref{sec:intro}, Chen et al.~\cite{Chen2000} have an impractical fully pipelined implementation. Coleman et al.~\cite{coleman2008processor} have add/sub taking 3 cycles to complete, but chose to duplicate rather than pipeline the unit, and mention that the latency is dominated by memory (ROM) access. Arnold~\cite{Arnold2003} provides a fully pipelined add/sub unit, but with a ``quick'' instruction version that allows the instruction to complete in either 4 or 6 cycles if it avoids the subtraction critical region. On the other hand, uniformity may increase latency, as different pipe stages are restricted to different ROM segments.

When combining an efficient LNS multiply with the penalty of addition for linear algebra, recent work by Popoff et al.~\cite{Popoff:2016:HLN:2971808.2972130} show an energy penalty of 1.84$\times$ over IEEE 754 float32 (using naive sum of add and mul energies), a 4.5$\times$ area penalty for the entire LNS ALU, and mention 25\% reduced performance for linear algebra kernels such as GEMM. Good LNS use cases likely remain workloads with high multiply-to-add ratios.

\subsection{ELMA/FLMA logarithmic arithmetic}\label{sec:flma}

The ELMA (\textit{exact log-linear multiply-add}) technique~\cite{jhj} is a logarithmic arithmetic that avoids Gaussian logarithms.
It was shown that an 8-bit ELMA implementation with extended dynamic range from posit-type encodings~\cite{Gustafson2017} is more energy efficient in 28 nm CMOS than 8/32-bit integer multiply-add (as used in neural network accelerators). It achieved similar accuracy as integer quantization on ResNet-50 CNN~\cite{he2016deep} inference on the ImageNet validation set~\cite{russakovsky2015imagenet}, simply with float32 parameters converted via round-to-nearest only and all arithmetic in the ELMA form. Significant energy efficiency gains over IEEE 754 float16 multiply-add were also shown, though much higher precision was then impractical.

We describe ELMA and its extension to FLMA (\textit{floating point log-linear multiply-add}).
In ELMA, mul/div/root/power is in log domain, while add/sub is in linear domain with fixed point arithmetic. Let $p(x')$ convert log domain $x'$ (with $E$ integer and $F$ fractional log bits) to linear domain, and $q(y)$ convert linear domain $y$ to log domain. $p(x')$ and $q(y)$ are both approximate conversions (LNS values are irrational). $p(x')$ produces fixed point (ELMA) or floating point (FLMA); in base-2 FLMA, we obtain $p(x') = \{ \pm 2^{\lfloor x' \rfloor} 2^{(x' - \lfloor x' \rfloor)}~\text{or}~ 0\}$, yielding a linear domain floating point exponent and significand. $p(x')$ can increase precision by $\alpha$ bits, with the exponential evaluated to $y_{\text{bits}} = f + \alpha$ fractional bits. Unique conversion for base-2 requires $\alpha \geq 1$, as the minimum derivative of $2^x$, $x \in [0, 1)$ is less than 1.

  FLMA approximates the linear domain sum $\sum_i x_i$ on the log domain $x'_i$ as $q(\sum_i p(x'_i))$. $q(\cdot)$ uses the floating point exponent as the log domain integer portion, and evaluates $\log_2$ on the significand, back to the required $F$ log domain fractional bits. The fixed or floating point accumulator can use a different fractional precision $A$ ($A \geq F + \alpha$) than $p(\cdot)$, in which case $q(\cdot)$ can consider $F + \beta$ linear domain MSB fractional bits of $A$ with rounding for the reverse conversion. $q(\cdot)$ is similarly unique only when $\beta \geq 1$. Typically we have $\alpha = \beta \geq 1$, and $A = F + \alpha$. As $\alpha, \beta$ increase, we converge to exact LNS add/sub. As with LNS, if add/sub is the only operation, ELMA/FLMA does not make sense. It is tailored for linear algebra sums-of-products; conversion errors are likely to be uncorrelated in use cases of interest (Sections~\ref{sec:add-sub} and~\ref{sec:multiple-sum}), but is substantially efficient over floating point at a 1:1 multiply-to-add ratio (Section~\ref{sec:syn}).


Unlike LNS, a ELMA design (and FLMA, depending upon floating point adder latency) can be easily pipelined and accept a new summand every cycle for accumulation without resource duplication (\textit{e.g.}, LNS ROMs). Furthermore, accumulator precision $A$ can be (much) greater than log domain $F$; in LNS this requires increasing Gaussian logarithm precision to the accumulator width. These properties make ELMA/FLMA excellent for inner product, where many sums of differing magnitudes may be accumulated.
FLMA is related to~\cite{Lai1991}, except that architecture is oriented around a linear domain floating point representation such that all mul/div/root is done with a log conversion to LNS, the log domain operation, and an exp conversion back to linear domain. Their log/exp conversions were further approximated with linear interpolation. Every mul/div/root operation thus included the error introduced by both conversions.

\subsection{Dual-base logarithmic arithmetic}\label{sec:dual}

ELMA/FLMA requires accurate calculation of the fractional portions of $p(x')$ and $q(y)$. Section~\ref{sec:exp} shows calculation of $e^x$ and $\ln(x)$ more accurately for the same resources versus $2^x$ and $\log_2(x)$. While Gaussian logarithms can be computed irrespective of base, FLMA requires an accessible base-2 exponent to carry over as a floating point exponent. A base-$e$ representation does not easily yield this. 

    An alternative is a variation of multiple base arithmetic by Dimitrov et al.~\cite{Dimitrov99}, allowing for more than one base $b_i$ (one of which is usually 2 and the others are any positive real number), with exponents $x_i$ as small integers, producing a representation $\pm \prod_i b^{x_i}_i$ (or zero). We instead use a representation $\pm 2^a e^b$ (or zero) with $a \in \mathbb{Z}$ (encoded in $E$ bits), $b \in [0, \ln(2))$ (encoded as an $F$ bit fixed point fraction).
      $e^b$ when evaluated yields a FLMA floating point significand in the range $[1, 2)$, which we will refer to as the \textit{Euler significand}. The product of any two of these values $2^a e^b \times 2^c e^d = 2^{a+c} e^{b+d}$ has $e^{b+d} \in [1, 4)$ and $(b + d) \in [0, 2\ln(2))$. For division, $e^{b-d} \in (0.5, 2), (b - d) \in (-\ln(2), \ln(2))$. We no longer have a unique representation when we do not limit the base-$e$ exponent to $[0, \ln(2))$; for example, $2^1 e^{0.4} \times 2^{-1} e^{0.3} = 2^0 e^{0.7} = 2^1 e^{(0.7 - \ln(2))}$.

              We call a base-$e$ exponent in the range $[0, \ln(2))$ a \textit{normalized Euler significand}. Normalization subtracts (or adds) $\ln(2)$ from the base-$e$ exponent and increments (or decrements) the base-2 exponent as necessary to obtain a normalized significand.
                There are two immediate downsides to this. First, we do not use the full encoding range; our base-$e$ exponent is encoded as a fixed point fraction, but we only use $\approx$ 69.3\% of the values. Encoding a precision/dynamic range tradeoff with the unused portion as in~\cite{Gustafson2017} could be considered. The second downside is considered in Section~\ref{sec:mul-div}.

\begin{figure}
    \vspace{-0.6cm}
  \begin{center}
    \includegraphics[width=1.0\linewidth]{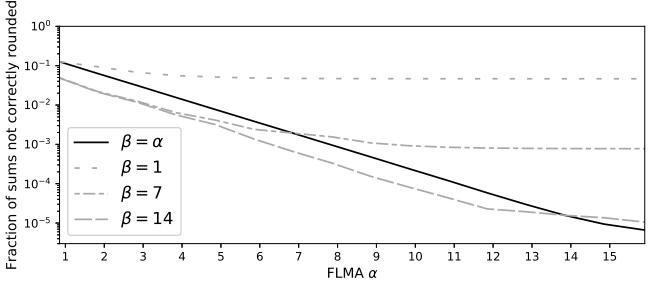}
  \end{center}
  \vspace{-0.5cm}
  \caption{\label{fig:add_ulp}
    FLMA sum error of log domain $x + y$ (with $x, y \in [1, 2)$) as a function of $\alpha, \beta$. All configurations have $\leq 1$ log ulp error, except for $\alpha = \beta = 1$ with $\leq 2$ log ulp error.
  }
\end{figure}

\section{FLMA analysis}\label{sec:flma_analysis}

We investigate dual-base $\pm 2^a e^b$ arithmetic with FLMA parameters $E=8$, $F=23$ (roughly IEEE 754 binary32 equivalent), exp $p=\ell=27+\alpha, r=2, I=13+\alpha$, log $p=\ell=27+\beta, I=14+\beta, r=3, s=9$, accumulator $A = F+\alpha$ for choice of $\alpha, \beta$. We call this log32 FLMA. For relative error, units in the last place in a fractional log domain representation we call \textit{log ulp}. For instance, 5 (base-$e$) log ulps are between $2^1 e^{\texttt{b0.0001}}$ and $2^1 e^{\texttt{b0.0110}}$, where \texttt{b0.0110} is the binary fixed point fraction 0.375.

\subsection{Multiply/divide accuracy}\label{sec:mul-div}

LNS and single base FLMA have 0 log ulp mul/div error, but dual-base FLMA can produce a non-normalized significand, requiring add/sub of a rounding of $\ln(2)$, introducing slight error ($\approx$ 0.016 log ulp for $F=23$).
The extended exp algorithm can avoid this for multiply-add with an additional iteration and integer bits for $L_n$ and $E_n$, as $2\ln(2) < 1.56$. We would still require additional normalization of $e^{b+d} \in [1, 4)$ to a floating point significand in the range $[1, 2)$. The dropped bit is kept by enlarging the accumulator, or is rounded away. Normalization is still required if more than two successive mul/div operations are performed.

\subsection{Add/subtract accuracy}\label{sec:add-sub}

Given $q(p(x') + p(y'))$ where both $x'$ and $y'$ are the same sign (\textit{i.e.}, not strict subtraction), the error is bounded by twice maximum $p(\cdot)$ error, plus maximum floating point addition error and maximum $q(\cdot)$ error. In practice the worst case error is hard to determine without exhaustive search. Limiting ourselves to values in a limited range of $[1, 2)$, we evaluate log domain FLMA addition of $x + y$ for all $x \in [1, 2)$ and a choice of 64 random $y \in [1, 2)$ versus $\alpha, \beta$ in Figure~\ref{fig:add_ulp} ($2^{29}$ unique pairs per configuration). All $\alpha > 1, \beta \geq 1$ have max log ulp error $\leq 1$, and $\alpha = \beta \leq 1$ has max log ulp error $\leq 2$. With increased $\alpha, \beta$ there are exponentially fewer incorrectly rounded sums but the table maker's dilemma is a limiting factor. At $\alpha = \beta = 16$, about 0.0005\% of these sums remain incorrectly rounded to max 0.5 log ulp.

\begin{figure}[b]
    \vspace{-0.3cm}
  \begin{center}
    \includegraphics[width=1.0\linewidth]{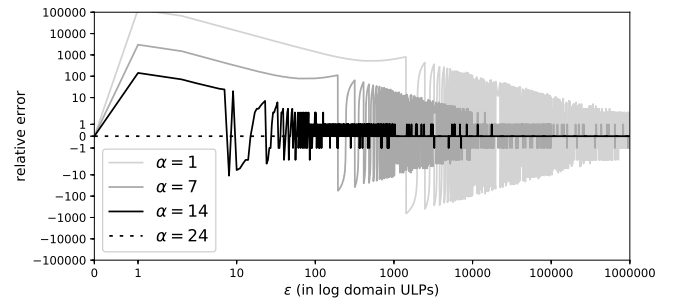}
  \end{center}
  \vspace{-0.5cm}
  \caption{\label{fig:catastrophe_rel}
    log32 FLMA catastrophic cancellation: relative (log ulp) error of $1 - (1 - \epsilon)$ in log domain, as a function of $\alpha$ ($\beta=1$ throughout).
  }
\end{figure}

For subtraction, catastrophic cancellation (a motivation for LNS co-transformation) still realizes itself. As with LNS, there is also a means of correction. While the issue appears with pairs of values very close in magnitude, consider linear domain $x = 1, y = 1 - \epsilon$, and evaluate $x - y$ with FLMA subtraction:
\begin{align*}
  x' &= +2^0 \times e^{\texttt{b0.00000000000000000000000}} \\
  y' &= +2^{-1} \times e^{\texttt{b0.10110001011100100001011}}
\end{align*}
The base-$e$ exponent of $y'$ here is 1 ulp below $\ln(2)$ rounded to $F=23$, and is thus our next lowest representable value from $1$. With FLMA subtraction at $\alpha = 1$:
\begin{align*}
  p(x') &= +2^0 \times \texttt{b1.0000000000000000000000} \\
  p(y') &= +2^{-1} \times \texttt{b1.1111111111111111111111} \\
  p(x') - p(y') &= +2^{-23} \times \texttt{b1.0000000000000000000000} \\
  &\approx 1.1920929 \times 10^{-7}
  \end{align*}
  Then back to log domain at $\beta = 1$:
\begin{align*}
  q(p(x') - p(y')) &= +2^{-23} \times e^{\texttt{b0.00000000000000000000000}}\\
  &\approx 1.1920929 \times 10^{-7}
  \end{align*}
If the calculation were done to $\leq$ 0.5 log ulp error, we get:
\begin{align*}
  q(p'(x') - p'(y')) &= +2^{-24} \times e^{\texttt{b0.10101101010100101000101}} \\
  & \approx 1.1730463 \times 10^{-7}
  \end{align*}
or an absolute error between the two of $\approx 1.905 \times 10^{-9}$, but off by 135,111 log ulp (distance from the $F=23$ bit rounding of $\ln(2)$).
In floating point, the rounded result would have error $\leq$ 0.5 ulp. However, as (almost) all of our log domain values have a linear domain infinite fractional expansion, in near cancellation with a limited number of bits, FLMA misses the extended expansion of the subtraction residual.

If reducing relative error is a concern, we can increase $\alpha$ for the log-to-linear conversion. This will provide more of the linear domain infinite fractional expansion, reducing relative error to $\leq 0.5$ log ulp almost everywhere if necessary (Figure~\ref{fig:catastrophe_rel}). Absolute error remains bounded throughout the cancellation regime, from $<10^{-8}$ at $\alpha=1$ to $<10^{-10}$ at $\alpha=14$.
We are not increasing the log precision $F$, but increasing the distinction in the linear domain between $1$ and $1 - \epsilon$. $q(y)$ can maintain a reduced $\beta$, with any remainder $A - (F + \beta)$ accumulator bits rounded off.

\begin{figure}[b]
  \begin{center}
    \includegraphics[width=1.0\linewidth]{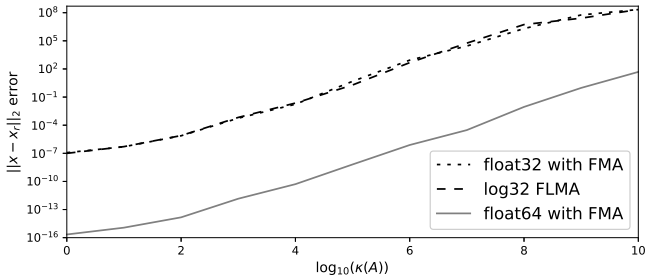}
  \end{center}
  \vspace{-0.5cm}
  \caption{\label{fig:mul-add-accuracy}
    Least squares $Ax = b$ solution average error via Householder QR, 64 $\times$ 64 matrices with condition number $\kappa(A)$ from 1 to $10^{10}$.
  }
\end{figure}


\subsection{Accuracy test via least squares QR solution}\label{sec:multiple-sum}

For a quick end-to-end accuracy test (encompassing sums, sums of products, multiplication, division and square root), we consider a least squares solution of $x$ in $Ax = b$ given various reference vectors $b$ and ill-conditioned reference matrices $A$. We control the condition number $\kappa(A)$ by generating random symmetric 64 $\times$ 64 matrices $M$ with entries $\sim U(-2, 2)$, taking a SVD decomposition $M = U \Sigma V$, scaling the largest to smallest singular value ratio by $\kappa$ to produce $\Sigma_\kappa$, then producing a reference $A = U \Sigma_\kappa V$ and $b$ ($b_i \sim U(0, 1)$) in float64 arithmetic. $A, b$ are rounded to $A', b'$ and QR factorization $A' = Q' R'$ is performed via the Householder algorithm in the arithmetic under analysis, and $x$ is recovered by backsolving $R' x = {Q'}^T b'$~\cite{golub13}. A reference $x_r$ from float64 $A, b$ is similarly calculated via MATLAB VPA to 64 decimal digits, and we present the error $\|x - x_r\|_2$. Figure~\ref{fig:mul-add-accuracy} shows this error across a sweep of condition number $\kappa(A)$ by powers of 10, with 5 trials for each condition number. Note that log32 FLMA remains roughly even against float32 (sometimes superior, sometimes inferior) despite the approximate nature of the design, even at high condition number $\kappa(A)=10^{10}$.


\section{Arithmetic synthesis}\label{sec:syn}

We compare 7 nm area, latency and energy against IEEE 754 floating point without subnormal support. A throughput of $T$ refers to a module accepting a new operation every $T$ clock cycles ($T=1$ is fully pipelined), while latency of $L$ is cycles to first result or pipeline length. Table~\ref{tab:arith-syn} shows basic arithmetic operations with FLMA parameters the same as Section~\ref{sec:flma_analysis} with $\alpha = \beta = 1$.
Note that the general LNS pattern of multiply energy being significantly lower but add/sub significantly higher still holds. Add/sub are two-operand, so this implementation includes two $p(\cdot)$ and one $q(\cdot)$ converters, and none will be actively gated in a fully utilized pipeline (they are all constantly switching). Naive sum of multiply with add energy lead to higher results as compared to floating point. However, as mentioned earlier, it is easier to efficiently pipeline FLMA add/sub compared to LNS add/sub.

  The situation changes when we consider a multiply-accumulate, perhaps the most important primitive for linear algebra. Table~\ref{tab:mul-add-syn} shows FLMA modules for 128-dim vector inner product with a $T=1$ inner loop, comparing against floating point FMA. The float64 comparison is against FLMA $E=11$, $F=52$, exp/log $\alpha = \beta = 1$, $p = \ell = 59$, $I = 29$, exp $r=2$, log $r=3, s = 9$, accumulator $A = 53$, called log64 FLMA.
  The benefit of the FLMA design can be seen in this case; log domain multiplication, $p(\cdot)$ conversion and floating point add is much lower energy than a floating point FMA. As with LNS or FLMA addition, a single multiply-add with a log domain result would be inefficient, but in running sum cases (multiply-accumulate), the $q(\cdot)$ overhead is deferred and amortized over all work, and this conversion (unlike the inner loop) need not be fully pipelined.
  Using a combinational MCP for this $q(\cdot)$ with data gating when inactive saves power and area, at the computational cost of 2 additional cycles for throughput. Increased accumulator precision ($A$ independent of $\alpha, \beta$) is also possible at minimal computational cost, as this only affects the floating point adder.

\begin{table}[t]
{
  \caption{\label{tab:arith-syn}
    Fully pipelined ($T=1$) arithmetic synthesis
}}
\centering \begin{tabular}{|l|l|l|l|}

\hline
\textbf{Type} & \textbf{Latency} & \textbf{Area} $\mu \text{m}^2$ & \textbf{Energy/op} pJ \\
\hline

float32 add/sub & 1 & 138.4 & 0.274\\
\textbf{log32 FLMA add/sub} & \textbf{7} & \textbf{1577 (11.4$\times$)} & \textbf{1.768 (6.45$\times$)} \\
\hline
float32 mul & 1 & 248.4 & 0.802\\
\textbf{log32 FLMA mul} & \textbf{1} & \textbf{40.2 (0.16$\times$)} & \textbf{0.080 (0.10$\times$)} \\
\hline
float32 FMA & 1 & 481.2 & 1.443\\
\textbf{log32 mul-add core} & & & \\ 
$p(x' + y') + a$, no $q(\cdot)$  & \textbf{3} & \textbf{706.5 (1.47$\times$)} & \textbf{0.586 (0.41$\times$)} \\
\hline
\end{tabular}

\end{table}

\begin{table}[t]
  \vspace{-0.3cm}
{
  \caption{\label{tab:mul-add-syn}
    $N=128$ inner product multiply-add synthesis results
}}
\centering \begin{tabular}{|l|l|l|l|}

\hline
\textbf{Type} & \textbf{Throughput} & \textbf{Area} $\mu \text{m}^2$ & \textbf{Energy/op} pJ \\

\hline
float32 FMA & 130 & 591.0 & 1.542 \\
\textbf{log32 FLMA} & & & \\
$q(\sum^{128}_{i=1} p(x'_i + y'_i))$ & \textbf{135} & \textbf{1271 (2.15$\times$)} & \textbf{0.668 (0.43$\times$)} \\
\hline
float64 FMA & 131 & 1787.3 & 5.032 \\
\textbf{log64 FLMA} & & & \\
$q(\sum^{128}_{i=1} p(x'_i + y'_i))$ & \textbf{144} & \textbf{6651 (3.72$\times$)} & \textbf{1.104 (0.22$\times$)} \\
\hline
\end{tabular}

\end{table}



\section{Conclusion}

\begin{figure}[t]
  \begin{center}
    \includegraphics[width=0.60\linewidth]{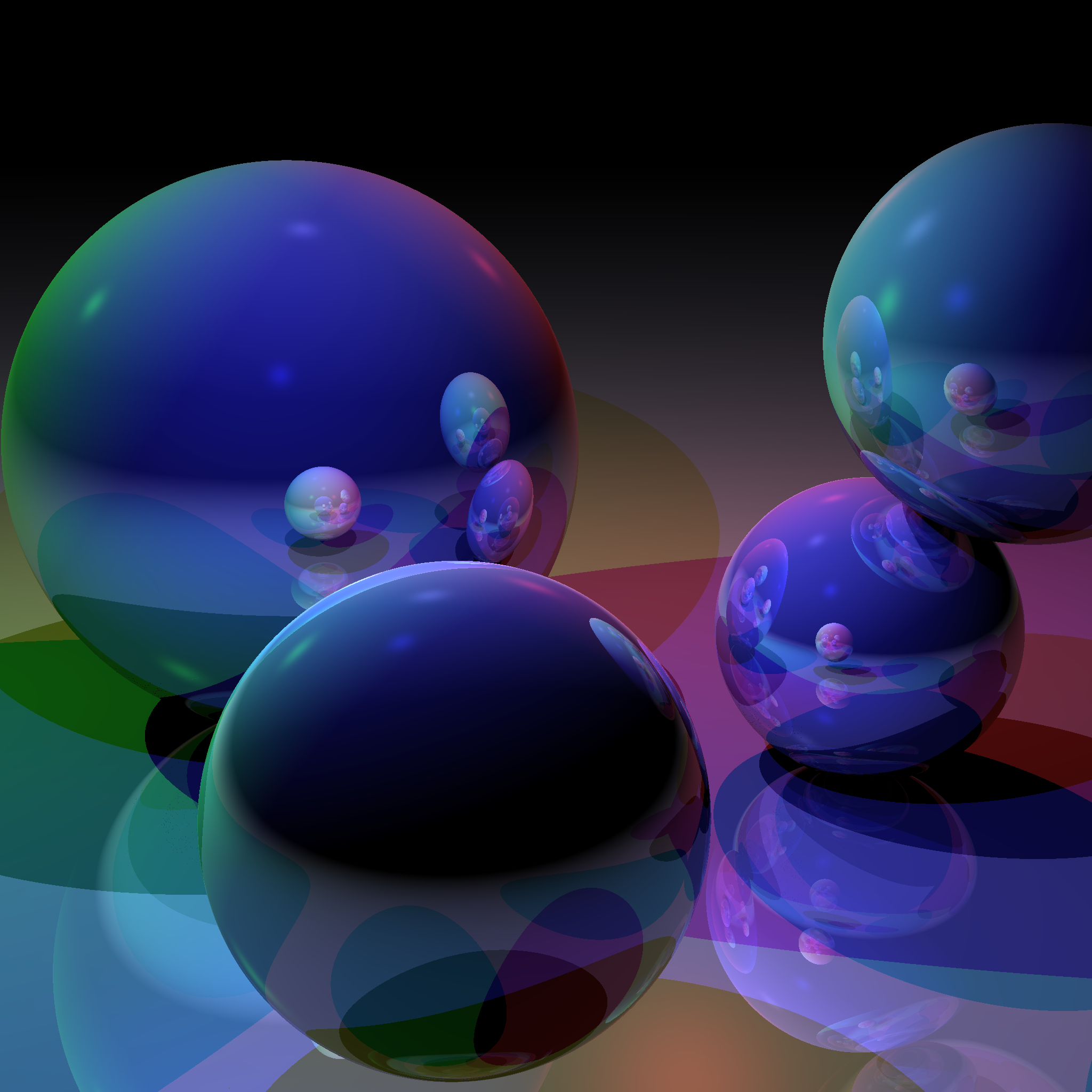}
  \end{center}
  \vspace{-0.5cm}
  \caption{\label{fig:raytrace}
    2048$\times$2048 raytracing done entirely in dual base FLMA arithmetic (Section~\ref{sec:flma_analysis} parameters with $\alpha = \beta = 1$). Pixel $\{ \pm 2^a e^b~\text{or}~0\}$ values clamped and rounded to nearest even integers for RGB output.
  }
\end{figure}

Modern applications of computer vision, graphics (Figure~\ref{fig:raytrace}) and machine learning often need energy efficient high precision arithmetic. We present an novel dual-base logarithmic arithmetic applicable for linear algebra kernels found in these applications. This is built on efficient implementations of $e^x$ and $\ln(x)$, useful in their own right, leveraging numerical integration with truncated mul/div.
While the arithmetic is approximate and without strong relative error guarantees unlike LNS or floating point arithmetic, it is extendible to arbitrary precision, easily pipelinable and retains moderate to low relative error and low absolute error. The area/power tradeoff is certainly not appropriate for many designs, but can provide a useful alternative to high precision floating or fixed point arithmetic when aggressive quantization is impractical.

\textsc{Acknowledgments} We thank Synopsys for their permission to publish results on our research obtained by using their tools with a popular 7 nm semiconductor technology node.

{
\bibliographystyle{IEEETran}
\bibliography{./IEEEabrv,./egbib}
}

\end{document}